\documentclass[12pt,oneside,leqno]{amsproc}
\usepackage[utf8]{inputenc}
\usepackage[russian]{babel}
\usepackage{amssymb,amsmath,amsthm}
\renewcommand{\subsection}{\refstepcounter{subsection}%
\par\bigskip\noindent\textbf{\upshape\thesubsection. }}
\renewcommand{\subsubsection}{\refstepcounter{subsubsection}%
\par\medskip\noindent\textbf{\upshape\thesubsubsection.  }}
\renewcommand{\paragraph}{\refstepcounter{paragraph}%
\par\smallskip\noindent\textbf{\upshape\theparagraph. }}

\numberwithin{equation}{subsection}

\renewcommand{\thesubsection}{\arabic{subsection}}
\renewcommand{\thesubsubsection}{\arabic{subsection}.\arabic{subsubsection}}

\newcommand{\Wo}{{\raisebox{0.2ex}{\(\stackrel{\circ}{W}\)}}{}}
\newcommand{\Ho}{\mathfrak H_0}
\newcommand{\ind}{\mathop{\operatorname{ind}}}
\newcommand{\im}{\operatorname{im}}

\newcommand{\supp}{\operatorname{supp}}

\makeatletter
\renewcommand{\@makefnmark}{\hbox{\small\mathsurround=0cm%
${}\hspace{0.04cm}{}^{\@thefnmark)}$}}
\renewcommand{\@makefntext}[1]{\parindent=1em\noindent\hbox to 1.8em{%
\hss${}^{\@thefnmark)}$}\,#1}
\makeatother

\title{Асимптотика собственных значений задачи высшего чётного порядка
с дискретным самоподобным весом}
\author{А.~А.~Владимиров, И.~А.~Шейпак\footnote{Работа поддержана РФФИ,
грант \No~10-01-00423.}}
\begin{document}
\renewcommand{\proofname}{{\upshape Д\,о\,к\,а\,з\,а\,т\,е\,л\,ь\,с\,т\,в\,о.}}
\begin{abstract}
В статье изучается вопрос об асимптотике спектра граничной задачи
\begin{gather*}
    (-1)^n\,y^{(2n)}-\lambda\rho y=0,\\
    y^{(k)}(0)=y^{(k)}(1)=0,\qquad 0\leqslant k<n,
\end{gather*}
в случае, когда порядок \(2n\) уравнения удовлетворяет неравенству \(n>1\),
а вес \(\rho\in\Wo_2^{-1}[0,1]\) представляет собой обобщённую производную
самоподобной функции \(P\in L_2[0,1]\) нулевого спектрального порядка.
\end{abstract}
\begin{flushleft}
\normalsize УДК~517.984
\end{flushleft}
\maketitle
\markboth{}{}

\section{Введение}\label{par:1}
\subsection\label{pt:1}
Целью настоящей статьи является распространение результатов работы \cite{VSh3}
о спектральных асимптотиках граничной задачи
\begin{gather*}
    -y''-\lambda\rho y=0,\\
    y(0)=y(1)=0,
\end{gather*}
где вес \(\rho\in\Wo_2^{-1}[0,1]\) представляет собой обобщённую производную
самоподобной функции \(P\in L_2[0,1]\) нулевого спектрального порядка, на случай
граничной задачи
\begin{gather}\label{eq:1}
    (-1)^n\,y^{(2n)}-\lambda\rho y=0,\\ \label{eq:2}
    y^{(k)}(0)=y^{(k)}(1)=0,\qquad 0\leqslant k<n,
\end{gather}
отвечающей тому же классу весовых функций при \(n>1\). Используемая в настоящей
статье техника в основном совпадает с развитой в работе \cite{VSh3}. Это, однако,
не означает, что переход от частного случая \(n=1\) к общему носит чисто
механический характер. Дополнительные трудности, возникающие в задаче высших
порядков (в отличие от случая задачи Штурма--Лиувилля), будут указаны нами далее.

На протяжении всей статьи мы резервируем символ \(n\) для обозначения половины
порядка уравнения \eqref{eq:1}.

\subsection
Структура статьи имеет следующий вид. В \ref{par:2} приводятся необходимые
для дальнейшего сведения о самоподобных функциях нулевого спектрального порядка.
В \ref{par:3} даётся операторная трактовка задачи \ref{pt:1}\,\eqref{eq:1},
\ref{pt:1}\,\eqref{eq:2} и доказываются некоторые вспомогательные утверждения.
В \ref{par:4} устанавливаются основные результаты об асимптотике спектра задачи
\ref{pt:1}\,\eqref{eq:1}, \ref{pt:1}\,\eqref{eq:2}. Наконец, в \ref{par:5}
обсуждается конструктивный характер полученных результатов, а также приводятся
иллюстрирующие эти результаты итоги численных экспериментов.

%%%%%%%%%%%%%%%%%%%%%%%%%%%%%%%%%%%%%%% СФ %%%%%%%%%%%%%%%%%%%%%%%%%%%%%%%%%%%%%%%%

\section{Квадратично суммируемые самоподобные функции нулевого спектрального
порядка}\label{par:2}
\subsection\label{pt:2.1}
Пусть зафиксировано натуральное число \(N>1\), и пусть вещественные числа \(a_k>0\),
\(\beta_k\) и \(d_k\), где \(k=1,\ldots N\), удовлетворяют равенству
\[
	\sum\limits_{k=1}^N a_k=1.
\]
Указанному набору чисел можно поставить в соответствие непрерывное отображение
\(G:L_2[0,1]\to L_2[0,1]\) вида
\begin{equation}\label{eq:auxto}
	G(f)\rightleftharpoons\sum\limits_{k=1}^N\left\{\beta_k\cdot
	\chi_{(\alpha_{k-1},\alpha_k)}+d_k\cdot G_k(f)\right\},
\end{equation}
где использована следующая символика:
\begin{enumerate}
\item Через \(\alpha_k\) обозначены числа \(\alpha_0\rightleftharpoons 0\)
и \(\alpha_k\rightleftharpoons\alpha_{k-1}+a_k\) при \(k=1,\ldots N\).
\item Через \(\chi_{\Delta}\) обозначена характеристическая функция интервала
\(\Delta\), рассматриваемая как элемент пространства \(L_2[0,1]\).
\item Через \(G_k\), где \(k=1,\ldots N\), обозначены непрерывные линейные
отображения в пространстве \(L_2[0,1]\), действующие на произвольно фиксированную
функцию \(f\in L_2[0,1]\) согласно правилу
\[
	[G_k(f)](x)\rightleftharpoons\left\{\begin{array}{ll}
		f((x-\alpha_{k-1})/a_k)&\text{при }x\in (\alpha_{k-1},\alpha_k),\\
		0&\text{иначе.}
	\end{array}\right.
\]
\end{enumerate}
Отображения вида~\eqref{eq:auxto} будут далее называться \emph{операторами
подобия}. Имеют место следующие два простых факта:

\subsubsection\label{lem2:1}
{\itshape Действующий в пространстве \(L_2[0,1]\) оператор подобия \(G\) является
сжимающим в том и только том случае, когда выполняется неравенствo
\begin{equation}\label{eq:szim}
	\sum\limits_{k=1}^N a_k\,|d_k|^2<1.
\end{equation}
}

\subsubsection\label{sek2:1}
{\itshape Если выполняется неравенство~\eqref{eq:szim}, то решение \(f\in L_2[0,1]\)
уравнения \(G(f)=f\) существует и единственно.
}

\medskip
Функции, удовлетворяющие уравнению \(G(f)=f\) с некоторым сжимающим оператором
подобия \(G\), мы будем называть \emph{аффинно самоподобными}, или просто
самоподобными. Определяющие оператор \(G\) числа \(a_k\), \(\beta_k\) и \(d_k\),
где \(k=1,\ldots N\), мы будем при этом называть \emph{параметрами самоподобия}
функции \(f\).

\subsection
Нетривиальная\footnote{То есть не являющаяся кусочно-постоянной с конечным числом
точек разрыва.} самоподобная функция называется функцией \emph{нулевого спектрального
порядка}, если её параметры самоподобия удовлетворяют следующим двум условиям:
\begin{enumerate}
\item Среди чисел \(\beta_k\), где \(k=1,\ldots N\), по меньшей мере одно отлично
от нуля.
\item Среди чисел \(d_k\), где \(k=1,\ldots N\), в точности одно отлично от нуля.
\end{enumerate}

В дальнейших рассуждениях о самоподобных функциях нулевого спектрального порядка через
\(m\) мы будем обозначать натуральное число \(m\in [1,N]\), удовлетворяющее соотношению
\(d_m\neq 0\). Неравенство \ref{pt:2.1}\,\eqref{eq:szim} при этом превращается
в неравенство \(a_m|d_m|^2<1\), из выполнения которого с очевидностью следует также
выполнение неравенства \(a_m^{2n-1}|d_m|<1\).

\subsection
Более подробные св\'{е}дения о квадратично суммируемых самоподобных функциях
могут быть найдены в работах \cite{Sh}, \cite{Sh2}.

%%%%%%%%%%%%%%%%%%%%%%%%%%%%%%%%%%%%%%%% ОТЗ %%%%%%%%%%%%%%%%%%%%%%%%%%%%%%%%%%%%%%%%%%%

\section{Операторная трактовка задачи и некоторые вспомогательные
утверждения}\label{par:3}
\subsection\label{pt:3.1}
Через \(\mathfrak H\) мы далее будем обозначать пространство Соболева
\(\Wo_2^n[0,1]\), снабжённое скалярным произведением
\[
    \langle y,z\rangle\rightleftharpoons\int\limits_0^1 y^{(n)}
	\overline{z^{(n)}}\,dx.
\]
Через \(\mathfrak H'\) мы при этом будем обозначать пространство, двойственное
к \(\mathfrak H\) относительно \(L_2[0,1]\), то есть получаемое пополнением
пространства \(L_2[0,1]\) по норме
\[
    \|y\|_{\mathfrak H'}\rightleftharpoons\sup\limits_{\|z\|_{\mathfrak H}=1}
    \left|\int\limits_0^1 y\overline{z}\,dx\right|.
\]
Непосредственно из определения пространства \(\mathfrak H'\) вытекает возможность
непрерывного продолжения сопряжённого к оператору вложения \(J:\mathfrak H\to
L_2[0,1]\) оператора \(J^*:L_2[0,1]\to\mathfrak H\) до изометрии
\(J^+:\mathfrak H'\to\mathfrak H\).

Аналогично использованной в работе \cite{VSh3} трактовке задачи Штурма--Лиувилля,
в качестве операторной модели задачи \ref{par:1}.\ref{pt:1}\,\eqref{eq:1},
\ref{par:1}.\ref{pt:1}\,\eqref{eq:2} мы будем рассматривать линейный пучок
\(T_{\rho}:\mathbb C\to\mathcal B(\mathfrak H,\mathfrak H')\) ограниченных
операторов, удовлетворяющий тождеству
\begin{equation}\label{eq:3.1}
    (\forall\lambda\in\mathbb C)\,(\forall y\in\mathfrak H)\qquad
    \langle J^+T_{\rho}(\lambda)y,y\rangle=\int\limits_0^1\left(|y^{(n)}|^2+
    \lambda P\cdot(|y|^2)'\right)\,dx.
\end{equation}
Через \(P\) здесь, как и ранее, обозначена квадратично суммируемая обобщённая
первообразная весовой функции \(\rho\in\Wo_2^{-1}[0,1]\).

\subsection
Через \(\ind D\) мы далее будем обозначать отрицательный индекс инерции действующего
в некотором гильбертовом пространстве \(\mathfrak E\) ограниченного эрмитова оператора
\(D\), то есть точную верхнюю грань размерностей подпространств \(\mathfrak M\subseteq
\mathfrak E\), удовлетворяющих условию
\[
    (\exists\varepsilon>0)\,(\forall y\in\mathfrak M)\qquad
    \langle Dy,y\rangle_{\mathfrak E}\leqslant-\varepsilon\,\|y\|^2_{\mathfrak E}.
\]

\subsection\label{pt:3:h}
Через \(\Ho\) мы далее будем обозначать замыкание линейной оболочки системы
собственных функций пучка \(T_{\rho}\). Имеют место следующие три факта:

\subsubsection\label{prop:h1}
{\itshape Ортогональным дополнением \(\mathfrak H\ominus\Ho\) подпространства
\(\Ho\) является множество функций \(y\in\mathfrak H\), обращающихся в нуль
на носителе \(\supp\rho\subseteq [0,1]\) весовой функции \(\rho\).
}

\begin{proof}
Из определения \ref{pt:3.1}\,\eqref{eq:3.1} и общей теории самосопряжённых
операторов в гильбертовом пространстве следует, что искомым ортогональным
дополнением является множество функций \(y\in\mathfrak H\), удовлетворяющих условию
\[
	(\forall z\in\mathfrak H)\qquad
	\int\limits_0^1 P\cdot (y\overline z)'\,dx=0.
\]
Данное условие равносильно указанному в формулировке доказываемого утверждения.
\end{proof}

\subsubsection\label{prop:h3}
{\itshape Пусть \(\lambda\) "--- вещественное число. Тогда существует
подпространство \(\mathfrak M\subseteq\Ho\) размерности
\(\ind J^+T_{\rho}(\lambda)\), удовлетворяющее условию
\[
	(\exists\varepsilon>0)\,(\forall y\in\mathfrak M)\qquad
	\langle J^+T_{\rho}(\lambda)y,y\rangle\leqslant
	-\varepsilon\,\|y\|^2_{\mathfrak H}.
\]
}

Для доказательства утверждения \ref{prop:h3} достаточно заметить, что для любых
функции \(y\in\mathfrak H\) и её ортогональной проекции \(z\in\mathfrak H_0\)
выполняется неравенство
\[
	\langle J^+T_{\rho}(\lambda)z,z\rangle\leqslant
	\langle J^+T_{\rho}(\lambda)y,y\rangle.
\]

\subsubsection\label{prop:h2}
{\itshape Для любой изолированной точки \(\xi\in\supp\rho\) существует
и единственна функция \(\varphi_{\xi}\in\Ho\), удовлетворяющая равенству
\(\varphi_{\xi}(\xi)=1\) и обращающаяся в нуль на множестве
\(\supp\rho\setminus\{\xi\}\).
}

\begin{proof}
Зафиксируем произвольную функцию \(f\in\mathfrak H\), удовлетворяющую равенству
\(f(\xi)=1\) и обращающуюся в нуль на множестве \(\supp\rho\setminus\{\xi\}\).
На роль искомой функции \(\varphi_{\xi}\) может быть теперь выбрана ортогональная
проекция функции \(f\) на подпространство \(\Ho\) [\ref{prop:h1}]. Для завершения
доказательства остаётся лишь заметить, что любая функция \(y\in\Ho\) однозначно
определяется своим ограничением на множество \(\supp\rho\) [\ref{prop:h1}].
\end{proof}

\subsection\label{pt:3.4}
Введём в рассмотрение два подпространства \(\mathfrak H_1\subseteq\Ho\)
и \(\mathfrak H_2\subseteq\Ho\), определяемые следующим образом:
\begin{enumerate}
\item Подпространство \(\mathfrak H_1\) образовано всевозможными функциями
\(y\in\Ho\), обращающимися в нуль на множестве \(\supp\rho\setminus(\alpha_{m-1},
\alpha_m)\).
\item Подпространство \(\mathfrak H_2\) представляет собой линейную оболочку
функций \(\varphi_{\xi}\) [\ref{prop:h2}], отвечающих всевозможным точкам
\(\xi\in\supp\rho\setminus (\alpha_{m-1},\alpha_m)\).
\end{enumerate}
Рассмотрим два линейных пучка \(A:\mathbb C\to\mathcal B(\mathfrak H_1,
\mathfrak H_1)\) и \(C:\mathbb C\to\mathcal B(\mathfrak H_2,\mathfrak H_2)\)
ограниченных операторов, удовлетворяющие тождествам
\begin{align}\label{eq:3.2}
	(\forall\lambda\in\mathbb C)\,(\forall y\in\mathfrak H_1)\qquad
	\langle A(\lambda)y,y\rangle&=\int\limits_0^1\left(|y^{(n)}|^2+
	\lambda P\cdot (|y|^2)'\right)\,dx,\\ \notag
	(\forall\lambda\in\mathbb C)\,(\forall y\in\mathfrak H_2)\qquad
	\langle C(\lambda)y,y\rangle&=\int\limits_0^1\left(|y^{(n)}|^2+
	\lambda P\cdot (|y|^2)'\right)\,dx,
\end{align}
а также оператор \(B:\mathfrak H_1\to\mathfrak H_2\), удовлетворяющий тождеству
\begin{equation}\label{eq:3.4}
    (\forall y\in\mathfrak H_1)\,(\forall z\in\mathfrak H_2)\qquad
    \langle By,z\rangle=\int\limits_0^1 y^{(n)}\overline{z^{(n)}}\,dx.
\end{equation}
Имеют место следующие два факта:

\subsubsection\label{st3:3}
{\itshape Существует неотрицательный оператор \(E:\mathfrak H_1\to\mathfrak H_1\)
конечного ранга, для которого независимо от выбора значения \(\lambda>0\)
выполняется равенство
\[
    \ind [A(\lambda)+E]=\ind J^+T_{\rho}(a_m^{2n-1}d_m\,\lambda).
\]
}

\begin{proof}
Рассмотрим оператор \(S:\mathfrak H\to\mathfrak H\) вида
\[
	[Sy](x)\rightleftharpoons\left\{\begin{array}{ll}
		y((x-\alpha_{m-1})/a_m)&\text{при }
		x\in(\alpha_{m-1},\alpha_m),\\
		0&\text{иначе,}
	\end{array}\right.
\]
а также оператор \(R:\mathfrak H_1\to\mathfrak H\), сопоставляющий каждой функции
\(y\in\mathfrak H_1\) ортогональную проекцию на подпространство \(\Ho\)
функции вида
\[
	z(x)\rightleftharpoons\psi(x)y(\alpha_{m-1}+a_mx),
\]
где \(\psi\in W_2^n[0,1]\) "--- произвольно фиксированная функция, тождественно
равная \(1\) на некоторой окрестности множества \(\supp\rho\) и обращающаяся
в нуль на некоторой окрестности множества \(\{0,1\}\setminus\supp\rho\).
Из утверждения \ref{prop:h1} и факта самоподобия функции \(P\) следует, что
для любой функции \(y\in\mathfrak H_1\) справедливы следующие два положения:
\begin{enumerate}
\item Ортогональная проекция функции \(SRy\) на подпространство \(\Ho\) совпадает
с функцией \(y\).
\item Функция \(Ry\) есть единственный элемент подпространства \(\Ho\),
ортогональная проекция \mbox{\(S\)-об}\-ра\-за которого на это же подпространство
совпадает с функцией \(y\).
\end{enumerate}
Соответственно, оператор \(E\rightleftharpoons R^*S^*SR-1\) является
неотрицательным. Из утверждения \ref{prop:h1} и факта самоподобия функции \(P\)
также следует, что на имеющем конечную коразмерность подпространстве функций
\(y\in\mathfrak H_1\), тождественно равных нулю вне интервала \((\alpha_{m-1},
\alpha_m)\), выполняется равенство \(SRy=y\). Тем самым, оператор \(E\) имеет
конечный ранг.

Наконец, из определений \ref{pt:3.1}\,\eqref{eq:3.1} и \eqref{eq:3.2} с учётом
утверждения \ref{prop:h1} и факта самоподобия функции \(P\) легко выводится
справедливость тождества
\[
	(\forall y\in\mathfrak H_1)\qquad \langle [A(\lambda)+E]y,y\rangle=
	a_m^{1-2n}\cdot\langle J^+T_{\rho}(a_m^{2n-1}d_m\,\lambda)Ry,Ry\rangle.
\]
Объединяя это тождество с утверждением \ref{prop:h3} и вытекающим из сказанного
ранее равенством \(\im R=\Ho\), убеждаемся в справедливости доказываемого
утверждения.
\end{proof}

\subsubsection\label{st3:4}
{\itshape Пусть \(\lambda\) "--- вещественное число, не принадлежащее спектру пучка
\(C\). Тогда выполняется равенство
\[
    \ind J^+T_{\rho}(\lambda)=\ind [A(\lambda)-B^*C^{-1}(\lambda)B]+
    \ind C(\lambda).
\]
}

\begin{proof}
Прямым вычислением легко устанавливается, что для любых функций
\(y\in\mathfrak H_1\) и \(z\in\mathfrak H_2\) выполняется равенство
\[
    \langle J^+T_{\rho}(\lambda)\,(y+z),(y+z)\rangle=
    \langle [A(\lambda)-B^*C^{-1}(\lambda)B]y,y\rangle+
    \langle C(\lambda)u,u\rangle,
\]
где положено \(u\rightleftharpoons z+C^{-1}(\lambda)By\). Объединяя это равенство
с утверждением \ref{prop:h3} и соотношением \(\Ho=\mathfrak H_1\dotplus
\mathfrak H_2\), убеждаемся в справедливости доказываемого утверждения.
\end{proof}

\subsection
Рассмотрим величины \(\zeta_k\), где \(k=1,\ldots N-1\), имеющие вид
\[
    \zeta_k\rightleftharpoons\left\{
        \begin{aligned}
            &\beta_m-\beta_{m-1}+d_m\beta_1&&\text{при }k=m-1,\\
            &\beta_{m+1}-\beta_m-d_m\beta_N&&\text{при }k=m,\\
            &\beta_{k+1}-\beta_k&&\text{иначе.}
        \end{aligned}
    \right.
\]
Обозначим также через \(\mathrm Z_{\pm}\) две величины
\[
    \mathrm Z_{\pm}\rightleftharpoons\#\{k\in [1,N-1]:\pm\zeta_k>0\}.
\]
Имеют место следующие два факта:

\subsubsection\label{4:1}
{\itshape Для любого достаточно большого вещественного числа \(\lambda>0\) выполняется
равенство
\[
    \ind C(\lambda)=\mathrm Z_+.
\]
}

Справедливость утверждения \ref{4:1} немедленно вытекает из легко проверяемого
тождества
\begin{equation}\label{eq:3}
    (\forall\lambda\in\mathbb R)\,(\forall y\in\mathfrak H_2)\qquad
    \langle C(\lambda)y,y\rangle=\|y\|^2_{\mathfrak H}-
    \lambda\sum\limits_{k=1}^{N-1}\zeta_k\cdot|y(\alpha_k)|^2.
\end{equation}

\subsubsection\label{4:2}
{\itshape Пусть выполнено равенство \(\mathrm Z_++\mathrm Z_-=N-1\). Тогда
для любого достаточно большого вещественного числа \(\lambda>0\) оператор
\(C(\lambda)\) является ограниченно обратимым, причём при \(\lambda\to+\infty\)
справедлива асимптотика
\[
    \|C^{-1}(\lambda)\|=O(\lambda^{-1}).
\]
}

Справедливость утверждения \ref{4:2} также представляет собой несложное следствие
тождества \eqref{eq:3}.

\subsection
Имеют место следующие два факта:

\subsubsection\label{prop:3:7:0}
{\itshape Пусть \(\mathfrak E\) "--- конечномерное гильбертово пространство,
\(D:\mathfrak E\to\mathfrak E\) "--- неотрицательный оператор,
а \(F:\mathfrak E\to\mathfrak E\) "--- эрмитов оператор. Пусть также
\(\{\mu_k\}_{k=1}^r\) и \(\{\lambda_k\}_{k=1}^r\), где положено
\(r\rightleftharpoons\operatorname{rank} F\) "--- списки сосчитанных
с учётом кратности собственных значений пучков \(1-\lambda F\) и \(1+D-
\lambda F\), соответственно. Тогда выполняется неравенство
\[
	\prod\limits_{k=1}^r\dfrac{\lambda_k}{\mu_k}\leqslant\det(1+D).
\]
}

\begin{proof}
Разложим пространство \(\mathfrak E\) в прямую сумму \(\ker F\oplus\im F\)
и рассмотрим отвечающие этому разложению блочно-матричные представления
\[
	D=\begin{pmatrix}D_{11}&D_{12}\\ D_{21}&D_{22}\end{pmatrix},\qquad
	F=\begin{pmatrix}0&0\\ 0&F_{22}\end{pmatrix}.
\]
Заметим, что спектр пучка \(1-\lambda F\) совпадает со спектром пучка
\(1-\lambda F_{22}\), а спектр пучка \(1+D-\lambda F\) совпадает со спектром
пучка \(1+D_{22}-D_{21}(1+D_{11})^{-1}D_{12}-\lambda F_{22}\). Соответственно,
имеют место равенства
\begin{align*}
	\prod\limits_{k=1}^r\mu_k&=\dfrac{1}{\det F_{22}},\\
	\prod\limits_{k=1}^r\lambda_k&=\dfrac{\det[1+D_{22}-D_{21}
	(1+D_{11})^{-1}D_{12}]}{\det F_{22}},
\end{align*}
а потому и равенства
\[
	\prod\limits_{k=1}^r\dfrac{\lambda_k}{\mu_k}=\det[1+D_{22}-
	D_{21}(1+D_{11})^{-1}D_{12}]=\dfrac{\det(1+D)}{\det(1+D_{11})}.
\]
Отсюда и из факта неотрицательности операторов \(D\) и \(D_{11}\) немедленно
вытекает справедливость доказываемого утверждения.
\end{proof}

\subsubsection\label{prop:3:7:1}
{\itshape Пусть \(\mathfrak E\) "--- сепарабельное гильбертово пространство,
\(D:\mathfrak E\to\mathfrak E\) "--- неотрицательный оператор конечного ранга,
а \(F:\mathfrak E\to\mathfrak E\) "--- вполне непрерывный эрмитов оператор.
Пусть также \(\{\mu_k\}_{k=1}^{\infty}\) и \(\{\lambda_k\}_{k=1}^{\infty}\) "---
последовательности занумерованных в порядке возрастания (с учётом кратности)
положительных собственных значений пучков \(1-\lambda F\) и \(1+D-\lambda F\),
соответственно. Тогда последовательность частичных произведений бесконечного
произведения
\begin{equation}\label{eq:prod}
	\prod\limits_{k=1}^{\infty}\dfrac{\lambda_k}{\mu_k}
\end{equation}
является неубывающей и ограниченной.
}

\begin{proof}
Заметим, что при любом \(\lambda\in\mathbb R\) выполняется неравенство
\(\ind [1+D-\lambda F]\leqslant\ind [1-\lambda F]\). Заметим также, что все
положительные собственные значения пучков \(1-\lambda F\) и \(1+D-\lambda F\)
имеют отрицательный тип (см., например, \cite{LSY}, \cite{Vl:2003}). Из общих
вариационных принципов для самосопряжённых оператор-функций (см. там же) потому
следует, что при любом \(k\geqslant 1\) выполняется соотношение \(\lambda_k
\geqslant\mu_k\). Тем самым, последовательность частичных произведений
бесконечного произведения \eqref{eq:prod} является неубывающей.

Далее, зафиксируем последовательность \(\{Q_l\}_{l=1}^{\infty}\) имеющих конечный
ранг ортопроекторов, сходящуюся в смысле сильной операторной топологии
к единичному оператору и удовлетворяющую тождеству \(Q_lDQ_l\equiv D\). Обозначим
через \(\{\mu_{k,l}\}_{k=1}^{\infty}\) и \(\{\lambda_{k,l}\}_{k=1}^{\infty}\)
последовательности\footnote{Частичные, то есть не предполагающие определённости
своих членов при каждом значении индекса \(k\geqslant 1\).} занумерованных
в порядке возрастания (с учётом кратности) положительных собственных значений
пучков \(1-\lambda Q_lFQ_l\) и \(1+D-\lambda Q_lFQ_l\), соответственно.
Аналогичным образом, через \(\{\mu_{-k,l}\}_{k=1}^{\infty}\)
и \(\{\lambda_{-k,l}\}_{k=1}^{\infty}\) обозначим последовательности занумерованных
в порядке убывания отрицательных собственных значений тех же пучков.
Из упоминавшихся выше вариационных принципов вытекает, что количество \(r_{+,l}\)
положительных и количество \(r_{-,l}\) отрицательных собственных значений пучка
\(1-\lambda Q_lFQ_l\) совпадают с таковыми для пучка \(1+D-\lambda Q_lFQ_l\),
причём выполняются неравенства
\begin{align*}
	&&&(\forall k\in [1,r_{+,l}])&\lambda_{k,l}&\geqslant\mu_{k,l},&&\\
	&&&(\forall k\in [1,r_{-,l}])&\lambda_{-k,l}&\leqslant\mu_{-k,l}.&&
\end{align*}
С учётом этого обстоятельства соотношение
\begin{flalign*}
	&&\prod\limits_{k=1}^{r_{+,l}}\dfrac{\lambda_{k,l}}{\mu_{k,l}}\times
	\prod\limits_{k=1}^{r_{-,l}}\dfrac{\lambda_{-k,l}}{\mu_{-k,l}}&\leqslant
	\det(1+D)&\text{[\ref{prop:3:7:0}]}&
\end{flalign*}
означает, что при любом выборе индекса \(r\leqslant r_{+,l}\) выполняются
неравенства
\[
	\prod\limits_{k=1}^{r}\dfrac{\lambda_{k,l}}{\mu_{k,l}}\leqslant\det(1+D).
\]
При помощи предельного перехода отсюда немедленно выводится справедливость
не зависящих от выбора индекса \(r\geqslant 1\) оценок
\[
	\prod\limits_{k=1}^r\dfrac{\lambda_k}{\mu_k}\leqslant\det(1+D),
\]
означающих ограниченность последовательности частичных произведений бесконечного
произведения \eqref{eq:prod}.
\end{proof}

%%%%%%%%%%%%%%%%%%%%%%%%%%%%%%%%%%%%%%%% ОР %%%%%%%%%%%%%%%%%%%%%%%%%%%%%%%%%%%%%%%%

\section{Основные результаты}\label{par:4}
\subsection
Имеют место следующие три факта:

\subsubsection\label{tm:1}
{\itshape Пусть выполняются соотношения \(d_m>0\), \(\mathrm Z_+>0\)
и \(\mathrm Z_++\mathrm Z_-=N-1\). Тогда существуют вещественные числа \(\tau_l>0\),
где \(l=1,\ldots\mathrm Z_{+}\), для которых последовательность
\(\{\lambda_k\}_{k=1}^{\infty}\) занумерованных в порядке возрастания (с учётом
кратности) положительных собственных значений задачи \ref{par:1}.\ref{pt:1}\,%
\eqref{eq:1}, \ref{par:1}.\ref{pt:1}\,\eqref{eq:2} удовлетворяет при \(k\to\infty\)
асимптотическому соотношению
\[
    \lambda_{l+k\mathrm Z_{+}}=\tau_l\cdot (a_m^{2n-1}d_m)^{-k}\cdot (1+o(1)).
\]
}

\begin{proof}
Согласно утверждению \ref{par:3}.\ref{4:2}, найдётся вещественное число
\(\lambda_0>0\), для которого при любом \(\lambda>\lambda_0\) будет выполняться
неравенство \(\|C^{-1}(\lambda)\|<\lambda_0/(3\lambda)\). Отсюда и из очевидным
образом получамой на основе определения \ref{par:3}.\ref{pt:3.4}\,\eqref{eq:3.4}
оценки \(\|B\|\leqslant 1\) следует, что при любом \(\lambda>\lambda_0\) будут
выполняться неравенства
\begin{equation}\label{eq:4}
    \ind [A(\lambda)+\lambda_0/(3\lambda)]\leqslant
    \ind [A(\lambda)-B^*C^{-1}(\lambda)B]\leqslant
    \ind [A(\lambda)-\lambda_0/(3\lambda)].
\end{equation}
С использованием немедленно вытекающих из определения
\ref{par:3}.\ref{pt:3.4}\,\eqref{eq:3.2} равенств
\[
    A(\lambda\pm\lambda_0/2)=\left(1\pm\lambda_0/(2\lambda)\right)\cdot
    \left[A(\lambda)\mp\dfrac{\lambda_0}{2\lambda\pm\lambda_0}\right]
\]
из оценок \eqref{eq:4} легко выводятся оценки
\begin{equation}\label{eq:5}
    \ind A(\lambda-\lambda_0/2)\leqslant
    \ind [A(\lambda)-B^*C^{-1}(\lambda)B]\leqslant
    \ind A(\lambda+\lambda_0/2).
\end{equation}

Обозначим теперь через \(\{\mu_k\}_{k=1}^{\infty}\) последовательность
занумерованных в порядке возрастания (с учётом кратности) положительных
собственных значений пучка \(A\). Из оценок \eqref{eq:5}, утверждений
\ref{par:3}.\ref{st3:4}, \ref{par:3}.\ref{4:1} и вариационных принципов
для самосопряжённых оператор-функций (см., например, \cite{LSY}, \cite{Vl:2003})
следует выполнение при всех \(k\gg 1\) неравенств
\[
	\mu_k-\lambda_0/2\leqslant\lambda_{k+\mathrm Z_+}\leqslant
	\mu_k+\lambda_0/2,
\]
а потому и асимптотического соотношения
\begin{equation}\label{eq:10}
	\dfrac{\lambda_{k+\mathrm Z_+}}{\mu_k}=1+O(\lambda_{k+\mathrm Z_+}^{-1}).
\end{equation}
С другой стороны, согласно утверждениям \ref{par:3}.\ref{st3:3}
и \ref{par:3}.\ref{prop:3:7:1}, последовательность частичных произведений
бесконечного произведения
\begin{equation}\label{eq:11}
	\prod\limits_{k=1}^{\infty}\dfrac{(a_m^{2n-1}d_m)\cdot\mu_k}{\lambda_k}
\end{equation}
является монотонной и ограниченной. Комбинируя этот факт с асимптотикой
\eqref{eq:10}, устанавливаем, что независимо от выбора параметров \(l\geqslant 1\)
и \(\varepsilon>0\) справедлива асимптотика
\[
	\lambda_{l+k\mathrm Z_+}^{-1}=O\bigl([(1+\varepsilon)\cdot
	a_m^{2n-1}d_m]^k\bigr),
\]
означающая сходимость бесконечного произведения
\[
	\prod\limits_{k=1}^{\infty}\dfrac{\lambda_{l+k\mathrm Z_+}}{%
	\mu_{l+(k-1)\mathrm Z_+}}.
\]
Для завершения доказательства остаётся теперь лишь объединить последний факт
с уже упомянутым ранее фактом монотонности и ограниченности последовательности
частичных произведений бесконечного произведения \eqref{eq:11}.
\end{proof}

\subsubsection\label{tm:2}
{\itshape Пусть выполняются соотношения \(d_m>0\), \(\mathrm Z_->0\)
и \(\mathrm Z_++\mathrm Z_-=N-1\). Тогда существуют вещественные числа \(\tau_l>0\),
где \(l=1,\ldots\mathrm Z_{-}\), для которых последовательность
\(\{\lambda_{-k}\}_{k=1}^{\infty}\) занумерованных в порядке убывания (с учётом
кратности) отрицательных собственных значений задачи \ref{par:1}.\ref{pt:1}\,%
\eqref{eq:1}, \ref{par:1}.\ref{pt:1}\,\eqref{eq:2} удовлетворяет при \(k\to\infty\)
асимптотическому соотношению
\[
    \lambda_{-(l+k\mathrm Z_{-})}=-\tau_l\cdot (a_m^{2n-1}d_m)^{-k}\cdot (1+o(1)).
\]
}

Доказательство утверждения \ref{tm:2} проводится аналогичным доказательству
утверждения \ref{tm:1} способом.

\subsubsection\label{tm:3}
{\itshape Пусть выполняются соотношения \(d_m<0\) и \(\mathrm Z_++\mathrm Z_-=N-1\).
Тогда существуют вещественные числа \(\tau_l>0\), где
\(l=1,\ldots N-1\), для которых последовательность \(\{\lambda_k\}_{k=1}^{\infty}\)
занумерованных в порядке возрастания (с учётом кратности) положительных собственных
значений задачи \ref{par:1}.\ref{pt:1}\,\eqref{eq:1}, \ref{par:1}.\ref{pt:1}\,%
\eqref{eq:2} удовлетворяет при \(k\to\infty\) асимптотическому соотношению
\[
    \lambda_{l+k(N-1)}=\tau_l\cdot (a_m^{2n-1}\,|d_m|)^{-2k}\cdot (1+o(1)),
\]
а последовательность \(\{\lambda_{-k}\}_{k=1}^{\infty}\) занумерованных в порядке
убывания (с учётом кратности) отрицательных собственных значений задачи
\ref{par:1}.\ref{pt:1}\,\eqref{eq:1}, \ref{par:1}.\ref{pt:1}\,\eqref{eq:2}
удовлетворяет при \(k\to\infty\) асимптотическому соотношению
\[
    \lambda_{-(l+\mathrm Z_{-}+k(N-1))}=-\tau_l\cdot (a_m^{2n-1}\,|d_m|)^{-2k-1}\cdot
    (1+o(1)).
\]
}

Доказательство утверждения \ref{tm:3} также проводится аналогичным доказательству
утверждения \ref{tm:1} способом.

%%%%%%%%%%%%%%%%%%%%%%%%%%%%%%%%%%%%%% Примеры %%%%%%%%%%%%%%%%%%%%%%%%%%%%%%%%%%%%%

\section{Примеры и обсуждение}\label{par:5}
\subsection
\begin{table}[t]
\begin{center}
\begin{tabular}{|r|r|r@{\;}c@{\;}l|r@{\;}c@{\;}l|}
\hline
{\(l\)}&{\(k\)}&\multicolumn{3}{|c|}{\(\lambda_{l+2k}\)}&
\multicolumn{3}{|c|}{\(\lambda_{l+2k}/54^k\)}\\ \hline
1&0&\(2,86\cdot 10^2\)&\(\pm\)&\(1\%\)&286,10&\(\pm\)&\(10^{-2}\)\\
2&0&\(1,38\cdot 10^3\)&\(\pm\)&\(1\%\)&1377,99&\(\pm\)&\(10^{-2}\)\\
1&1&\(1,48\cdot 10^4\)&\(\pm\)&\(1\%\)&273,71&\(\pm\)&\(10^{-2}\)\\
2&1&\(6,83\cdot 10^4\)&\(\pm\)&\(1\%\)&1265,31&\(\pm\)&\(10^{-2}\)\\
1&2&\(7,91\cdot 10^5\)&\(\pm\)&\(1\%\)&271,33&\(\pm\)&\(10^{-2}\)\\
2&2&\(3,69\cdot 10^6\)&\(\pm\)&\(1\%\)&1264,04&\(\pm\)&\(10^{-2}\)\\
1&3&\(4,27\cdot 10^7\)&\(\pm\)&\(1\%\)&271,32&\(\pm\)&\(10^{-2}\)\\
2&3&\(1,99\cdot 10^8\)&\(\pm\)&\(1\%\)&1264,04&\(\pm\)&\(10^{-2}\)\\
\hline
\end{tabular}
\end{center}

\vspace{0.5cm}
\caption{Оценки первых собственных значений для случая \(n=2\), \(N=3\), \(a_1=a_2=
a_3=1/3\), \(m=3\), \(d_3=1/2\), \(\beta_1=0\), \(\beta_2=2/3\), \(\beta_3=1\).}
\label{tab:1}
\end{table}
В таблице \ref{tab:1} представлены результаты численных расчётов для первых восьми
положительных собственных значений задачи четвёртого порядка, весовой функцией
в которой выступает обобщённая производная квадратично суммируемой функции
с параметрами самоподобия \(N=3\), \(a_1=a_2=a_3=1/3\), \(m=3\), \(d_3=1/2\),
\(\beta_1=0\), \(\beta_2=2/3\), \(\beta_3=1\). В этом случае выполняются
равенства \(\zeta_1=2/3\), \(\zeta_2=1/3\), \(\mathrm Z_+=2\), \(\mathrm Z_-=0\).
Данные таблицы иллюстрируют утверждение \ref{par:4}.\ref{tm:1}.

\begin{table}[t]
\begin{center}
\begin{tabular}{|r|r|r@{\;}c@{\;}l|r@{\;}c@{\;}l|}
\hline
{\(l\)}&{\(k\)}&\multicolumn{3}{|c|}{\(-\lambda_{-(l+k)}\)}&
\multicolumn{3}{|c|}{\(-\lambda_{-(l+k)}/54^{k}\)}\\ \hline
1&0&\(3,70\cdot 10^2\)&\(\pm\)&\(1\%\)&369,75&\(\pm\)&\(10^{-2}\)\\
1&1&\(8,51\cdot 10^3\)&\(\pm\)&\(1\%\)&157,53&\(\pm\)&\(10^{-2}\)\\
1&2&\(4,58\cdot 10^5\)&\(\pm\)&\(1\%\)&157,20&\(\pm\)&\(10^{-2}\)\\
1&3&\(2,48\cdot 10^7\)&\(\pm\)&\(1\%\)&157,20&\(\pm\)&\(10^{-2}\)\\
\hline
\end{tabular}
\end{center}

\vspace{0.5cm}
\caption{Оценки первых собственных значений для случая \(n=2\), \(N=3\),
\(a_1=a_2=a_3=1/3\), \(m=3\), \(d_3=1/2\), \(\beta_1=\beta_3=0\),
\(\beta_2=-1\).}
\label{tab:2}
\end{table}
В таблице \ref{tab:2} представлены данные численных расчётов первых четырёх
отрицательных собственных значений задачи четвёртого порядка, весовой функцией
в которой выступает обобщённая производная квадратично суммируемой функции
с параметрами самоподобия \(N=3\), \(a_1=a_2=a_3=1/3\), \(m=3\), \(d_3=1/2\),
\(\beta_1=\beta_3=0\), \(\beta_2=-1\). В этом случае выполняются равенства
\(\zeta_1=-1\), \(\zeta_2=1\), \(\mathrm Z_+=\mathrm Z_-=1\). Данные таблицы
иллюстрируют утверждение \ref{par:4}.\ref{tm:2}.

\begin{table}[t]
\begin{center}
\begin{tabular}{|r|r|r@{\;}c@{\;}l|r@{\;}c@{\;}l|r@{\;}c@{\;}l|r@{\;}c@{\;}l|}
\hline
{\(l\)}&{\(k\)}&\multicolumn{3}{|c|}{\(\lambda_{l+2k}\)}&
\multicolumn{3}{|c|}{\(\lambda_{l+2k}/54^{2k}\)}&
\multicolumn{3}{|c|}{\(-\lambda_{-(l+1+2k)}\)}&
\multicolumn{3}{|c|}{\(-\lambda_{-(l+1+2k)}/54^{2k+1}\)}\\ \hline
1&0&\(3,04\cdot 10^2\)&\(\pm\)&\(1\%\)&304,08&\(\pm\)&\(10^{-2}\)&
\(1,61\cdot 10^4\)&\(\pm\)&\(1\%\)&299,04&\(\pm\)&\(10^{-2}\)\\
2&0&\(1,38\cdot 10^4\)&\(\pm\)&\(1\%\)&13820,11&\(\pm\)&\(10^{-2}\)&
\(7,43\cdot 10^5\)&\(\pm\)&\(1\%\)&13764,22&\(\pm\)&\(10^{-2}\)\\
1&1&\(8,72\cdot 10^5\)&\(\pm\)&\(1\%\)&299,00&\(\pm\)&\(10^{-2}\)&
\(4,71\cdot 10^7\)&\(\pm\)&\(1\%\)&299,00&\(\pm\)&\(10^{-2}\)\\
2&1&\(4,01\cdot 10^7\)&\(\pm\)&\(1\%\)&13764,02&\(\pm\)&\(10^{-2}\)&
\(2,17\cdot 10^9\)&\(\pm\)&\(1\%\)&13764,02&\(\pm\)&\(10^{-2}\)\\
1&2&\(2,54\cdot 10^9\)&\(\pm\)&\(1\%\)&299,00&\(\pm\)&\(10^{-2}\)&
\(1,37\cdot 10^{11}\)&\(\pm\)&\(1\%\)&299,00&\(\pm\)&\(10^{-2}\)\\
2&2&\(1,17\cdot 10^{11}\)&\(\pm\)&\(1\%\)&13764,02&\(\pm\)&\(10^{-2}\)&
\(6,32\cdot 10^{12}\)&\(\pm\)&\(1\%\)&13764,02&\(\pm\)&\(10^{-2}\)\\
\hline
\end{tabular}
\end{center}

\vspace{0.5cm}
\caption{Оценки первых собственных значений для случая \(n=2\), \(N=3\),
\(a_1=a_2=a_3=1/3\), \(m=3\), \(d_3=-1/2\), \(\beta_1=\beta_3=0\),
\(\beta_2=-1\).}
\label{tab:3}
\end{table}
В таблице \ref{tab:3} представлены данные численных расчётов первых шести
положительных и семи отрицательных (исключая первое) собственных значений задачи
четвёртого порядка, весовой функцией в которой выступает обобщённая производная
квадратично суммируемой функции с параметрами самоподобия \(N=3\),
\(a_1=a_2=a_3=1/3\), \(m=3\), \(d_3=-1/2\), \(\beta_1=\beta_3=0\), \(\beta_2=-1\).
В этом случае выполняются равенства \(\zeta_1=-1\), \(\zeta_2=1\),
\(\mathrm Z_+=\mathrm Z_-=1\). Данные таблицы иллюстрируют утверждение
\ref{par:4}.\ref{tm:3}.

Для получения вышеприведённого иллюстративного материала нами была
применена "--- в незначительным образом модифицированном виде "---
вычислительная методика, описанная в работе \cite{V}.

\subsection
Некоторые из рассуждений, проведённых нами в предыдущих параграфах, не являются
приемлемыми с точки зрения конструктивного направления в математике \cite{Sh:1962},
\cite{Ku:1973}. Это относится, в частности, к характеру установления утверждений
из пункта \ref{par:3}.\ref{pt:3:h}, связанному с привлечением теоремы
об ортогональном проектировании в гильбертовом пространстве. Несколько более
громоздкое рассуждение (на деталях которого мы здесь не останавливаемся)
позволило бы избежать опоры на следствия этой теоремы, не меняя существа дела.
Однако использование в ходе доказательств утверждений из~\ref{par:4} теоремы
Больцано--Вейерштрасса о сходимости монотонной ограниченной последовательности
вещественных чисел является в рамках применённого нами подхода, по-видимому,
неустранимым. Соответственно, коэффициенты \(\tau_l\) полученных асимптотик
оказываются с точки зрения конструктивного математического анализа так называемыми
\emph{псевдочислами}\footnote{Для которых эффективный метод построения сколь угодно
точных рациональных приближений либо неизвестен, либо прямо невозможен.}. Данное
обстоятельство коренным образом отличает результаты настоящей статьи от результатов
работы \cite{VSh3}, где соответствующие коэффициенты заведомо являлись
\emph{конструктивными вещественными числами}\footnote{Допускающими вычисление сколь
угодно точных рациональных приближений посредством единого алгорифма.}: в случае
задачи Штурма--Лиувилля оператор \(E\) из утверждения \ref{par:3}.\ref{st3:3}
является нулевым, что позволяет полностью избежать апелляции к утверждению
\ref{par:3}.\ref{prop:3:7:1}. Вопрос о возможности дальнейшего уточнения
конструктивного характера полученных асимптотик мы оставляем открытым.

%%%%%%%%%%%%%%%%%%%%%%%%%%%%%%%%%% Литература %%%%%%%%%%%%%%%%%%%%%%%%%%%%%%%%%%%%%

\end{document}